# On 2-digit and 3-digit Kaprekar's Routine

Zichen Wang, Wei Lu

**Abstract**

Kaprekar's Routine is an iteration process that, with each iteration, sorts the digits of a number in descending and ascending orders and uses their difference for the next iteration. In this paper, we successfully address the structure of the *fixed periodic sets* and the *maximum step of iteration* to enter a fixed set in three-digit case and two-digit case, with a focus on the latter. With three-digit, there exist a unique fixed set of cardinality of 1 or 2, depending on the parity of $m$, and the maximum step is determined by $m$. With two-digit, we offer some powerful properties of the fixed sets and a mechanical method to calculate the fixed sets. We also prove that the number of 2-factors in $m+1$ determines the maximum step (but has completely no influence on the structure of the fixed sets). This continues AtsushiYamagami's previous studies in 2018 (published on Journal of Number Theory) and completes it with a general solution.

# Contents

# 1 Introduction

## 1.1 Overview

In 1949, D. R. Kaprekar [1] first discovered the following process: takes any four-digit number with at least two different digits; arrange the digits in ascending and descending orders to generate two new numbers (not necessarily four-digit); subtract the smaller number from the bigger number; and iterate the process with the new number generated. Eventually, the process will generate 6174 within 7 iterations, and the number remains fixed as 7641-1467=6174.

$$3332-2333=999$$

$$9990-0999=8991$$

$$9981-1899=8082$$

$$8820-288=8532$$

$$8532-2358=6174$$

This process later becomes known as the Kaprekar's Routine, and 6174 is called a Kaprekar's Constant. Generalizing this process, we can see that Kaprekar's Routine holds with different digits and in different bases. For example, 495 is another Kaprekar's Constant with three digits in base ten.

Admittedly, Kaprekar's Routine began as a recreational exploration and now most commonly appears as "number black holes" in math magazines [6][7]. However, recent studies have revealed that two-digit Kaprekar's Routine is closely related to Mersenne primes [4], and the whole process has found independent applications to cryptographic encoding schemes [9].

Though there is growing attention on this subject, so far little research about Kaprekar's Routine has been conducted and limited papers are published. Many authors from previous studies still relied heavily on computer simulations and statistical methods to observe certain patterns of the Kaprekar's Routine [8], and a more theoretical explanation of the subject is in urgent need.

## 1.2 Acknowledgement

This work was done as an independent research project at Southeast University, Nanjing. The first author is grateful for the guidance from Dr. Lu, who is the first author's supervisor, the director of the mini-project, and the second author of this paper. He would also like to thank New York University for providing access to various data bases.

Contact information:
zichenatwork@163.com (Zichen Wang), luwei1010@seu.edu.cn (Wei Lu).

# 2 Discrete Dynamic System

## 2.1 Definitions

We first introduce the concept of *discrete dynamic system* [2][3], which consists of a set *S* and a self-map *f* on *S*. It helps us characterize the iteration process as a morphism *f* and the *m*-base *n*-digit numbers as a set *S*. For the sake of this paper, we will give directly a discrete dynamic system that describe specifically the Kaprekar's Routine. However, notice that the following definitions can be easily extended to describe general discrete dynamic systems.

Let *X* be a subset of natural numbers $\mathbb{N}$, *f* be self-map on *X*, that is

$$f : X \to X$$

$$x \to f(x)$$

**Definition 2.1.1.** Suppose $x \in X$, $t \in \mathbb{N}$, we define inductively

$$f^t(x) := \begin{cases} x & t = 0 \\ f(f^{t-1}(x)) & t > 0 \end{cases}$$

With the dynamic system set up, we can now define certain terms to our concern.

**Definition 2.1.2.** For a value $x \in X$, if there exists $s \in \mathbb{N}$, $t \in \mathbb{N}^+$, such that

$$f^{s+t}(x) = f^s(x) \qquad (2.1.1)$$

we say that *x* is *periodic* and *t* is a *period* of *x*.

**Definition 2.1.3.** For a value $x \in X$, we define

*the step S(x)* to be the smallest *s* that satisfies (2.1.1);

*the maximum step S(X)* to be the maximum *S(x)* for all $x \in X$;

*the minimal period T(x)* to be the smallest *t* that satisfies (2.1.1);

and *the fixed set* $K(x) := \{f^{S(x)+i}(x) \mid 0 \le i < T(x)\}$.

Intuitively, *step* describes the first time a value enter its *fixed set*, for $f^{S(x)}(x)$ is an element of *K(x)*. Also, the cardinality of the fixed set satisfies $|K(x)| = T(x)$.

## 2.2 Properties

With the definitions to help us set up the dynamic system, we can now give certain properties of the system.

**Lemma 2.2.1.** *If a value (after several iterations) enters its fixed set, its following iterations would only generate elements in the set.*

We will not proof this lemma here. However, it is worthy to point out that the *fixed set K(x)* nicely describes the iterations after $\{f^t(x)\}_{t\in\mathbb{N}}$ enters it (see the example below). Also, we can define a fixed set by *S(x)* and *T(x)* because the following relationship

$$f^{S(x)+T(x)}(x) = f^{S(x)}(x)$$

**Example 2.2.1.** In base 10, with 2 digits, *K*(15)={9,81,63,27,45}, while

$$15 \to 36 \to 27 \to 45 \to 9 \to 81 \to 63 \to 27 \to 45 \to 9 \to 81 \to 63 \to \ldots\ldots$$

**Lemma 2.2.2.** *For any base and any digit, all values x in the corresponding set X are periodic. In fact, all values will inevitably fall into and remain inside its fixed set.*

*Proof.* Suppose the numbers have *m*-base and *n*-digit, then the corresponding set *X* is a finite set with no more than $(m+1)^n$ elements. However, the number-generating process can iterate infinite time, and this guarantees that at least a number appears twice in the sequence and thus all numbers are *periodic*. For a specific number *x*, it will iterate *S(x)* times to first satisfy (2.1.1) and thus enter a period. Hence, $\{f^t(x)\}_{t\in\mathbb{N}}$ will inevitably enter *K(x)* and by Lemma 2.2.1., remain inside the fixed set. ∎

## 2.3 Questions

Given the definitions and properties of the dynamic system, we are thus able to ask the following questions. Let *X* be a set of *m*-base *n*-digit numbers, *f* a self-map on *X*,

1. For which *m*, *n*, there exist *fixed sets* with cardinality of 1?
   (i.e., the Kaprekar Constant)
2. Given a positive natural number *t*, for which *m*, *n*, there exist $x \in X$, *T(x)=t*?
3. For which *m*, *n*, $\forall x_1, x_2 \in X, T(x_1) = T(x_2)$ ?
   (i.e., there exist only one *minimal period*)
4. For which *m*, *n*, $\forall x_1, x_2 \in X, K(x_1) = K(x_2)$ ?
   (i.e., there exist only one *fixed set*)
5. Given *m*, *n*, can we tell *S(X)*?
6. Given *m*, *n*, and a value $x \in X$, can we tell *S(x)*?

# 3 Kaprekar's Routine

## 3.1 General Discussion

Consider the original four-digit Kaprekar's Routine. We can express any four-digit number as

$$\overline{abcd} = 1000a + 100b + 10c + d$$

Rearrange *a, b, c, d* in descending order and we get $a_1, a_2, a_3, a_4$. Thus, we can write one iteration as

$$f(\overline{abcd}) = \overline{a_1a_2a_3a_4} - \overline{a_4a_3a_2a_1}$$
$$= (1000a_1 + 100a_2 + 10a_3 + a_4) - (1000a_4 + 100a_3 + 10a_2 + a_1)$$
$$= 999a_1 + 90a_2 - 90a_3 - 999a_4$$

Notice that 0 always satisfies as a Kaprekar' Constant. We call this a *trivial Kaprekar's Constant*, and will ignore it in the following discussions. Notice also that the number becomes a multiply of 9 after one iteration and retains this property thereafter. This factor of 9 greatly limits the possibility of the final Kaprekar's Constant and is derived by $(10^k - 1)a_i$. Hence, we see that the base number greatly influences the final result. This serves as an intuition for the generalization of the Kaprekar's Routine in all bases.

## 3.2 Three-Digit Case

Surprisingly, the three-digit Kaprekar's Routine comes much simpler than the two-digit case, and thus we decide to discuss the three-digit case first. A three-digit number in base $m$ can be written as $\overline{abc}_m = m^2 a + mb + c$, and the set $X = \left\{ \sum_{i=1}^{3} a_i m^{i-1} \mid 0 \leq a_i < m \right\}$.

We continue the above discussion to figure out the common factor first.

**Lemma 3.2.1** *For any $x \in X$, $m^2-1$ divides $f(x)$. In fact,*

$$f(X) = \left\{ a(m^2 - 1) \mid 0 \leq a \leq m \right\}$$

*Proof.* Rearrange $x$'s digits and we get $a_0, a_1, a_2$ ($0 \leq a_i < m$) in descending order. Thus,

$$f(x) = \overline{a_1 a_2 a_3} - \overline{a_3 a_2 a_1}$$
$$= (m^2 a_1 + ma_2 + a_3) - (m^2 a_3 + ma_2 + a_1)$$
$$= (a_1 - a_3)(m^2 - 1)$$

Let $a = a_1 - a_3$. When $a_1, a_3$ change, $a$ can achieve any integer between 0 and $m$. ∎

At this point, we will shift our focuses on $a$ and supplement the following definition

**Definition 3.2.1.** Define $X_m := \{a \in \mathbb{N} \mid 0 \leq a \leq m\}$ and $f_3$ a self-map on $X_m$, that is

$$f_3 : X_m \to X_m$$
$$x \to f_3(x),$$

so that for any $a \in X_m$, there exist some $f_3(a) \in X_m$, $f(a(m^2 - 1)) = f_3(a)(m^2 - 1)$.

Now we only need to consider the period of the sequence $\{f_3^t(a)\}_{t \in \mathbb{N}}$, and at core the question is about solving the equation

$$f_3^t(a) = a, a \in X_m$$

**Lemma 3.2.2.** *When $a \neq 0$, we have* $f_3(a) = \begin{cases} a-1 & a \geq \dfrac{m+1}{2} \\ m-a & a < \dfrac{m+1}{2} \end{cases}$

*Proof.* We can write the m-base three-digit number as

$$a(m^2 - 1) = (a-1)m^2 + (m-1)m + (m-a) = \overline{(a-1)(m-1)(m-a)}_m$$

If $a - 1 \geq m - a$, the digits have order $m - 1 \geq a - 1 \geq m - a$, thus

$$f\left(a(m^2 - 1)\right) = ((m-1) - (m-a)) \cdot (m^2 - 1) = (a-1)(m^2 - 1),$$
$$f_3(a) = a - 1$$

If $a - 1 < m - a$,

$$f\left(a(m^2 - 1)\right) = ((m-1) - (a-1)) \cdot (m^2 - 1) = (m-a)(m^2 - 1),$$
$$f_3(a) = m - a \qquad \blacksquare$$

Observe the above result, we notice that for any $a$, $f_3(a)$ is bigger than $(m+1)/2$. Subsequent iterations decrease the number by one each time, and we can expect the number to stabilize at or oscillate around $m/2$. In addition, the constant decrease of one enables us to nicely predict the *step* of the number. Strictly verifying these observations, we arrive at the most important theorem of this section.

**Theorem 3.2.3.** *Suppose* $x \in X$, $f(x) = a(m^2 - 1)$ *for some* $a \in X_m$

(a) *If m is even, X has a unique non-trivial fixed set* $\{\dfrac{m}{2}(m^2 - 1)\}$, *with*

$$S(x) = \left| a - \dfrac{m+1}{2} \right| + \dfrac{1}{2} + \delta(x)$$

(b) *If m is odd, X has a unique non-trivial fixed set* $\{\dfrac{m-1}{2}(m^2 - 1), \dfrac{m+1}{2}(m^2 - 1)\}$

$$S(x) = \left| a - \dfrac{m+1}{2} \right| + \delta(x)$$

(c) $\delta(x) = \begin{cases} 1 & (m^2 - 1) \nmid x \\ 0 & (m^2 - 1) \mid x \end{cases}$

(d) $S(X) = \begin{cases} \dfrac{m}{2} + 1 & m \text{ is even} \\ \dfrac{m+1}{2} & m \text{ is odd} \end{cases}$

*Proof.*

If $a \geq \frac{m+1}{2}$, $f_3(a) = a - 1 \geq \frac{m-1}{2}$; if $a < \frac{m+1}{2}$, $f_3(a) = m - a > \frac{m-1}{2}$. Either way, $f_3(a) \geq \frac{m-1}{2}$.

When $f_3^t(a) > \frac{m+1}{2}$, $f_3^{t+1}(a) = f_3^t(a) - 1$ decreases monotonically until

(a) $f_3^t(a) = \frac{m}{2}$.

$f_3^t(a)$ enters its fixed set (point) as $f_3^{t+1}(a) = f_3^t(a) = \frac{m}{2}$.

Since $f_3^t(a)$ is always an integer, $m$ must be an even number.

(b) $f_3^t(a) = \frac{m-1}{2}$ or $\frac{m+1}{2}$.

$f_3^t(a)$ enters its fixed set as $f_3^{t+1}(a) = \frac{m+1}{2}$ or $\frac{m-1}{2}$

Since $f_3^t(a)$ is always an integer, $m$ must be an odd number.

At this point, the proof for $S(X)$ should be straightforward and is left for readers. Notice, however, that we need $\delta(x)$ to judge if $x$ is already a multiple of $(m^2-1)$. Otherwise, we need one iteration to first transfer $x$ into the form of $a(m^2-1)$. ∎

We conclude this section by offering two examples.

**Example 3.2.1.** In base 10, $X$ has a unique non-trivial fixed set {495}.

$S(X)=10/2+1=6$ establishes as $667 \to 99 \to 891 \to 792 \to 693 \to 594 \to 495 \to 495$.

**Example 3.2.2.** In base 13, $X$ has a unique non-trivial fixed set {(5,12,7), (7,12,5)}. $S(X)=(13+1)/2=7$ establishes as

$(1,0,1) \to (0,12,12) \to (11,12,1) \to (10,12,2) \to (9,12,3) \to (8,12,4) \to (7,12,5) \to (5,12,7)$.

## 3.3 Two-Digit Case: Structure of Fixed Sets

A two-digit number in base $m$ can be written as $\overline{ab}_m = ma + b$, and the set of all two-digit numbers $X = \{am + b \mid 0 \leq a, b < m\}$. Now, we want to continue our methodology with three-digit to factorize and thus simplify $f$ first.

**Lemma 3.3.1.** *For any $x \in X$, $m$-1 divides $f(x)$. In fact, $f(X) = \{a(m-1) \mid 0 \leq a \leq m\}$*

*Proof.* $f(x) = \overline{ab} - \overline{ba} = (ma + b) - (mb + a) = (m-1)(a - b)$ ($0 < b < a < m$)

($a$-$b$) can achieve any number between 0 and $m$. ∎

Similarly, we use $X_m$ defined in Definition 3.2.1. and supplement the definition of $f_2$.

**Definition 3.3.1.** Define $f_2$ a self-map on $X_m$, that is

$$f_2 : X_m \to X_m$$

$$x \to f_2(x),$$

so that for any $a \in X_m$, there exist some $f_2(a) \in X_m$, $f(a(m-1)) = f_2(a)(m-1)$.

Again, we shift our study to $f(a)$ and give the result of one iteration first.

**Lemma 3.3.2.** *When $a \neq 0$, we have $f_2(a) = |2a - m - 1|$.*

*Proof.* We can write the m-base two-digit number as

$$a(m-1) = (a-1)m + (m-a) = \overline{(a-1)(m-a)}_m$$

Thus, $f(a(m-1)) = |(a-1)-(m-a)|(m-1) = |2a-m-1|(m-1)$. ∎

This set the basis for the discussion of $f_2^t(a) = a, a \in X_m$.

**Lemma 3.3.3.** *X has a non-trivial fixed set of cardinality of 1 if and only if $3 | m+1$. In fact, such a fixed set is unique.*

*Proof.* To have a non-trivial fixed set of cardinality of 1, there must exist some $a \in X_m$, $f_2(a)=|2a-m-1|=a$. Solving this equation, we have $a=m+1$ or $a=(m+1)/3$. However, $a$ is smaller than $m$ and thus can only be $(m+1)/3$. ∎

**Lemma 3.3.4.** *X has a non-trivial fixed set of cardinality of 2 if and only if $5 | m+1$. In fact, such a fixed set is unique.*

*Proof.* To have a fixed set of cardinality of 2, there must exist some $a \in X_m$, $f_2^2(a) = a$.

$$f_2^2(a) = |2f_2(a) - m - 1| = a$$

$$f_2(a) = |2a - m - 1| = \frac{m+1}{2} \pm \frac{a}{2}$$

$$a = \frac{2 \pm 1}{2^2 \pm 1}(m+1) \qquad \text{(each } \pm \text{ is independent)}$$

$$a = \frac{3}{5}(m+1) \text{ or } \frac{1}{5}(m+1) \text{ or } (m+1) \text{ or } \frac{1}{3}(m+1)$$

$a$ is smaller than $m$ and thus cannot be $(m+1)$. Also, $a=(m+1)/3$ is a solution of $f(a)=a$, which is also a solution of $f_2^2(a) = a$ and thus appears here. Thus, we have only two values of $a$ remaining and they determine a unique fixed set of cardinality of 2. ∎

In general, suppose $f_2(a)=A$, then $a = \frac{1}{2}(m+1) \pm \frac{A}{2}$. Utilize repeatively this solution, we can solve the essential equation $f_2^t(a) = a$.

**Theorem 3.3.1.** *Given t, there exist $a \in X_m$, $f_2^t(a) = a$ if and only if $\gcd(m+1, 2^t \pm 1) > 1$. In fact, all non-zero solutions of $f_2^t(a) = a$ can be written as $a = \dfrac{\xi(m+1)}{2^t \pm 1}$.*

*($\xi$ can be any odd number smaller than $2^t$.)*

*Proof.*
$$a = f_2^t(a) = f_2(f_2^{t-1}(a))$$

$$f_2^{t-1}(a) = f_2(f_2^{t-2}(a)) = \frac{1}{2}(m+1) \pm \frac{a}{2}$$

$$f_2^{t-2}(a) = f_2(f_2^{t-3}(a)) = \frac{1}{2}(m+1) \pm \frac{1}{2}\left(\frac{1}{2}(m+1) + \frac{a}{2}\right)$$

$$= \left(\frac{1}{2} \pm \frac{1}{2^2}\right)(m+1) \pm \frac{a}{2^2}$$

Continuing this process and we will eventually get

$$a = \left(\frac{1}{2} \pm \frac{1}{2^2} \pm \frac{1}{2^4} \pm \ldots \pm \frac{1}{2^t}\right)(m+1) \pm \frac{a}{2^t}$$

$$a = \frac{2^{t-1} \pm 2^{t-2} \pm \ldots \pm 2 \pm 1}{2^t \pm 1}(m+1) \quad \text{(each $\pm$ is independent)}$$

$2^{t-1} \pm 2^{t-2} \pm \cdots 2 \pm 1$ can be any odd number between 0 and $2^t$. Also, there must be $\gcd(m+1, 2^t \pm 1) > 1$ for $a$ to be smaller than $m$. On the other hand, if $\gcd(m+1, 2^t \pm 1) > 1$, we can always let $\xi = (2^t \pm 1)/\gcd(m+1, 2^t \pm 1)$, and we can easily verify that $f_2^t(a) = a$ for such an $a$. ∎

This immediately explains Lemma 3.3.3 and Lemma 3.3.4. Since $2^1 + 1 = 3$ is a prime, $\gcd(m+1, 2^1 \pm 1) > 1$ implies that 3 is a factor of $m+1$. The case for $2^2 + 1 = 5$ is similar. Furthermore, the theorem offers us the discriminant of $f_2^t(a) = a$ and the property of the solutions. We can thus write down the general solutions of the equation.

**Theorem 3.3.2.** *Given t, the set of all solutions of $f_2^t(a) = a$ can be written as*

$$\left\{ a \in X_m \,\middle|\, a = \frac{m+1}{\gcd(m+1, 2^t+1)} \cdot \xi_1 \,\, (\xi_1 \text{ can be any odd number between 1 and } \gcd(m+1, 2^t+1)) \right\}$$

$$\cup \left\{ a \in X_m \,\middle|\, a = \frac{m+1}{\gcd(m+1, 2^t-1)} \cdot \xi_2 \,\, (\xi_2 \text{ can be any odd number between 1 and } \gcd(m+1, 2^t-1)) \right\}$$

*Proof.* By Theorem 3.3.1., all non-zero solutions of $f_2^t(a) = a$ can be written as

$$a = \frac{\xi(m+1)}{2^t \pm 1} = \frac{\xi\left((m+1)/\gcd(m+1, 2^t \pm 1)\right)}{(2^t \pm 1)/\gcd(m+1, 2^t \pm 1)}.$$

Since $\gcd\left((m+1)/\gcd(m+1, 2^t \pm 1), (2^t \pm 1)/\gcd(m+1, 2^t \pm 1)\right) = 1$, $\xi$ must be in form of

$$\frac{\xi'(2^t \pm 1)}{\gcd(m+1, 2^t \pm 1)}$$ for $a$ to be an integer. Thus, $a = \frac{\xi(m+1)}{2^t \pm 1} = \frac{(m+1)}{2^t \pm 1} \cdot \frac{\xi'(2^t \pm 1)}{\gcd(m+1, 2^t \pm 1)}$

$= \frac{\xi'(m+1)}{\gcd(m+1, 2^t \pm 1)}$. Notice that $a<m$, so $1 \leq \xi' < \gcd(m+1, 2^t \pm 1)$. ∎

With Theorem 3.3.2., we can write down all values of a certain *period*. However, from Lemma 3.3.3 and Lemma 3.3.4., we see that a solution for $f_2(a) = a$ is also a solution of $f_2^2(a) = a$. This leads us to further explore the *minimal period*.

Notice that the solutions of the equation are largely determined by $\gcd(m+1, 2^t \pm 1)$. If for a given $m$, there exists $t_1 < t_2$, $\gcd(m+1, 2^{t_1} \pm 1) = \gcd(m+1, 2^{t_2} \pm 1)$, then by Theorem 3.3.2., the solutions for $f_2^{t_1}(a) = a$ and $f_2^{t_2}(a) = a$ are exactly the same. This implies that $t_1 | t_2$, and while $t_2$ is a valid period, it cannot be a minimal period.

**Theorem 3.3.3.** *Given t, there exist $x \in X$, $T(x)=t$ if and only if there exist some $d|m+1$, $t$ is the smallest integer that satisfies $d|2^t+1$ or $d|2^t-1$.*

*Proof.* Sufficiency.

Suppose $\gcd(m+1, 2^t \pm 1) \geq d > 1$. By Theorem 3.3.1., $t$ is a period for some $x \in X$. If $T(x)=t'<t$, then again by Theorem 3.3.1., $\gcd(m+1, 2^{t'} \pm 1) > 1$. This contradicts that $t$ is the smallest integer that satisfies $d|2^t+1$ or $d|2^t-1$. Thus $T(x)=t$.

Necessity.

Let $d=\gcd(m+1, 2^t \pm 1)$. If for some $t'<t$, $d|2^{t'}+1$ or $d|2^{t'}-1$, then by Theorem 3.3.2., we can verify that all solutions for $f_2^t(a) = a$ are also solutions of $f_2^{t'}(a) = a$. This contradicts that $t$ is the minimal period. ∎

Here we list some interesting corollaries below (proof omitted). However, be aware that we can derive much more specific statements of the structure of the fixed sets based on Theorem 3.3.3. Examples include the condition for $X$ to have only one minimal period (but may have multiple fixed sets of the same period), the condition for a unique non-trivial fixed set, and the relationship between the cardinality of the fixed sets etc.

**Corollary 3.3.1.** *The number of 2-factors of $m+1$ does not affect the structure of the fixed sets (see Example 3.3.2. and Example 3.3.3.).*

**Corollary 3.3.2.** *If $m+1$ is a power of 2, then X has only the trivial fixed point 0.*

Now a few words on Theorem 3.3.3. The authors do not think we can push this conclusion any further, as we do not think it is possible to give an exact expression of the smallest $t$ in modern mathematics. However, Theorem 3.3.3 offers us a method to calculate all minimal periods for a given base, and by Theorem 3.3.2., we can write down precisely the values in each fixed set. Thus, we can calculate the fixed sets mechanically.

**Example 3.3.1.** In base 27, $m+1=28=2^2 \cdot 7$.

By Theorem 3.3.3., we examine the factors of 28 to get the only minimal period: 3.

By Theorem 3.3.2., we solve all elements in the fixed set: {4·26, 12·26, 20·26}.

(written in base 10 and in the form of $a \cdot (m\text{-}1)$)

**Example 3.3.2.** In base 14, $m+1=15=3 \cdot 5$.

By Theorem 3.3.3., we examine the factors of 15 to get the minimal periods: 1, 2, 4.

By Theorem 3.3.2., we solve all elements in the fixed sets: {5·13}; {3·13, 9·13}; {1·13, 7·13, 11·13, 13·13}.

**Example 3.3.3.** In base 59, $m+1=60=2^2 \cdot 3 \cdot 5$.

By Theorem 3.3.3., we examine the factors of 60 to get the minimal periods: 1, 2, 4.

By Theorem 3.3.2., we solve all elements in the fixed sets: {20·58}; {12·58, 36·58}; {4·58, 28·58, 44·58, 52·58}.

## 3.4 Two-Digit: Maximum Step

In last section, we studied in detail the structure of the fixed set. We will move on to study the *maximum step $S(X)$* to enter a fixed set in this section. Observe Example 3.3.1. and Example 3.3.2., we discover that all elements in the fixed sets have the same number of 2-factors with $m+1$. This offers us a crucial insight to the study of $S(X)$.

Let $p^r \| n$ denote $p^r | n$ but $p^{r+1} \nmid n$ (i.e., $n$ has $r$ 2-factors).

**Lemma 3.4.1.** *Suppose $2^r \| m+1$, then the union set of all fixed sets can be written as*

$$\left\{ a \in X_m \mid 2^r \| a \right\}$$

*Proof.* Re-statement: given $a \in X_m$, $a$ belongs to some fixed set if and only if $2^r \| a$.

Necessity. If $a$ belongs to some fixed set, then by Theorem 3.3.1, $a = \xi(m+1)/(2^t \pm 1)$. Notice that $\xi$ and $2^t \pm 1$ are both odd numbers, so $2^r \| m+1$ implies that $2^r \| a$.

Sufficiency. Suppose $m+1=2^r m'$ and $a=2^r a'$. Since 2 and $m'$ are coprime, by Euler's Theorem, $2^{\varphi(m')} \equiv 1 \pmod{m'}$ and thus $\gcd(2^{\varphi(m')}-1, m)=m'$. Now we can write

$$a = 2^r a' = \frac{m+1}{m'} \cdot a' = \frac{m+1}{\gcd(m+1, 2^{\varphi(m')}-1)} \cdot a'$$

By Theorem 3.3.2., $a$ belongs to some fixed set. ∎

In other words, we can say that the iteration process before entering a fixed set is the process of "gaining or losing 2-factors" until the number has precisely $r$ 2-factors. If there exist some nice pattern of 2-factors in the iteration process, then the problem would be solved.

**Theorem 3.4.1.** *Suppose $2^r \| m+1$, for any $x \in X$, $x$ will enter a fixed set within $r+2$ steps.*

*Proof.* Suppose after one iteration (at most), $f_2(x)$ is in form of $a(m-1)$, and $2^l \| a$.

Let $m=2^r \cdot m'$, $a=2^l \cdot a'$ ($m'$ and $a'$ are both odd numbers).

(a) If $r>l+1$,
$$f_2(a) = |2a-(m+1)| \text{ (by Lemma 3.3.1.)}$$
$$= |2^{l+1}a' - 2^r m'|$$
$$= 2^{l+1}|a' - 2^{r-l-1}m'|$$

$a'$ is odd and $2^{r-l-1}m'$ is even, so $|a' - 2^{r-l-1}m'|$ is odd.

Thus, $f_2(a)$ has one more 2-factor than $a$.

(b) If $r=l+1$,
$$f_2(a) = |2a-(m+1)|$$
$$= |2^{l+1}a' - 2^r m'|$$
$$= 2^r |a' - m'|$$

$a'$ is odd and $m'$ is odd, so $|a' - m'|$ is even.

Thus, $f_2(a)$ has more than $r$ 2-factors.

(c) If $r<l$,
$$f_2(a) = |2a-(m+1)|$$
$$= |2^{l+1}a' - 2^r m'|$$
$$= 2^r |2^{l-r-1}a' - m'|$$

$2^{l-r-1}a'$ is even and $m'$ is odd, so $|2^{l-r-1}a' - m'|$ is odd.

Thus, $f_2(a)$ has precisely $r$ 2-factors and is in some fixed set.

At most, $x$ needs one step to enter case (a), repeat case (a) $r-1$ times to enter case(b). Therefore we need at most $1+(r-1)+1+1=r+2$ steps. ∎

In fact, we can write down the expression of the step for a given x to enter a fixed set. The argument is almost identical to the proof above, and is thus omitted here.

**Example 3.4.1.** In base 14, $m+1=15=3\cdot5$.

$S(X)=0+2=2$ establishes as $21 \rightarrow 6\cdot13 \rightarrow 3\cdot13 \rightarrow 9\cdot13 \rightarrow 3\cdot13$.

**Example 3.4.2.** In base 59, $m+1=60=2^2\cdot3\cdot5$.

$S(X)=2+2=4$ establishes as $63 \rightarrow 3\cdot58 \rightarrow 54\cdot58 \rightarrow 48\cdot58 \rightarrow 36\cdot58 \rightarrow 12\cdot58 \rightarrow 36\cdot58$.

(Notice the change of number of 2-factors in the iteration process)

# 4 Concluding Remarks

From the original 6174 Kaprekar's Routine, we see that the common factor 9 for all numbers after one iteration greatly limits the possibility of the fixed set (point). This leads us to generalize our discussion to different bases, as well as our methodology of factoring out the common factor to focus on the remaining portion. In section 2, We introduced the concept of *discrete dynamic system* to mathematically describe the iteration process. This offers us a fundamental tool for our study of the Kaprekar's Routine.

In section 3, which is the main part of this paper, we focused our discussion on the structure of the fixed sets and the maximum step a number needs to enter a fixed set. With three-digit, we answered perfectly the above questions: there exist a unique fixed set of cardinality of 1 or 2, depending on the parity of *m*; the maximum step is also determined by *m*. With two-digit, the situation becomes more complex. While we do not think it is possible to write an expression of the fixed sets, we derive some important properties of the cardinality of fixed sets and the elements in them. Also, we offer a purely mechanical method to solve all the fixed sets for a given *m*. Lastly, we see that the property of $2^r \pm 1$ greatly determines the structure of the fixed sets, while the number of 2-factors of *m*+1 has no influence on the structure of the fixed sets, but determines the maximum step to enter a fixed set.

For future research, we suggest two potential directions here. One natural direction is to study Kaprekar's Routine with larger digits. As a general tool in the field of arithmetic dynamics, the discrete dynamic system is still applicable here. However, our method of factorization might not be as powerful as in lower digits since the property of *f(x)* weakens in higher digits. Another direction is the study of more specific structure of the fixed set in the two-digit case. Although a general expression might not be possible, we can still solve the expressions of the fixed sets for certain *m*. This involves more classical number theory, and might offer new insights to some unsolved problems (such as the Mersenne primes).

# 5 References

[1] D.R. Kaprekar. Another solitaire game. *Scripta Math*, 15:244-245, 1949.

[2] Joseph H. Silverman. *The Arithmetic of Dynamical Systems*, volume 241 of Graduate Texts in Mathematics. Springer, New York, 2007.

[3] Robert L. Benedetto and Joseph H. Silverman. Current Trends and Open Problems in Arithmetic Dynamics. *Bulletin of the American Mathematical Society*, 0273-0979(XX)0000-0.

[4] Atsushi Yamagami. On 2-adic Kaprekar constants and 2-digit Kaprekar distances. *Journal of Number Theory*, 185:257-280, 2018.

[5] Pat Devlin and Tony Zeng. Maximum distances in the four-digit Kaprekar process. arXiv: 2010.11756.

[6] Martin Gardner. Mathematical games. *Scientific American*, 232(3):112-117, 1975.

[7] Michael Ecker. Caution: Black Holes at Work. *New Scientist*, December 1992.

[8] Deniel Hanover. The Base Dependent Behavior of Kaprekar's Routine: A Theoretical and Computational Study Revealing New Regularities. arXiv: 1710.06308.

[9] Rishab G Nandan and Ritvika G Nandan. Multi-layer encryption employing Kaprekar routine and letter-proximity-based cryptograms, February 20 2020. US Patent App. 16/600,524.